\documentclass[11pt]{article}
\usepackage{graphicx}\usepackage{float}\usepackage{color}
\usepackage{amsmath}\usepackage{amsfonts}\usepackage{amssymb}
\providecommand{\U}[1]{\protect\rule{.1in}{.1in}}

\begin{document}

\begin{center}
{\Large \textbf{ Art of Unlocking }
\let\thefootnote\relax\footnotetext{\noindent
The work of the first author (O. O.) was supported by the grant PhyMath ANR-19-CE40-00. Part of the work of the second author (S.S.) has
been carried out at Skoltech and at IITP RAS. The support of Russian Science
Foundation  20-41-09009 is gratefully
acknowledged.}}

\vspace{.4cm} {\large \textbf{Oleg Ogievetsky$^{1,2,3}$
and Senya Shlosman$^{1,3,4}$}}

\vskip .3cm $^{1}$Aix Marseille Universit\'{e}, Universit\'{e} de
Toulon, CNRS, \\ CPT UMR 7332, 13288, Marseille, France

\vskip .05cm $^{2}$I.E.Tamm Department of Theoretical Physics,
Lebedev Physical Institute,
Leninsky prospect 53, 119991, Moscow, Russia

\vskip .05cm $^{3}$Inst. of the Information Transmission Problems, RAS,
Moscow, Russia

\vskip .05cm $^{4}$Skolkovo Institute of Science and Technology,
Moscow, Russia
\end{center}

\vskip .4cm
\begin{abstract}
We describe our recent results concerning the rigidity/unlockability properties of clusters of rigid bodies sliding over the unit sphere.
\end{abstract}

\section{Introduction: Balls}

The most known question in the area of movable/rigid configurations deals with
the FCC and HCP configurations of 12 unit balls touching the central unit ball
$\mathbb{B}$ (shown in red):

\begin{figure}[H]
\vspace{-.0cm}
\centering
\includegraphics[scale=0.068]{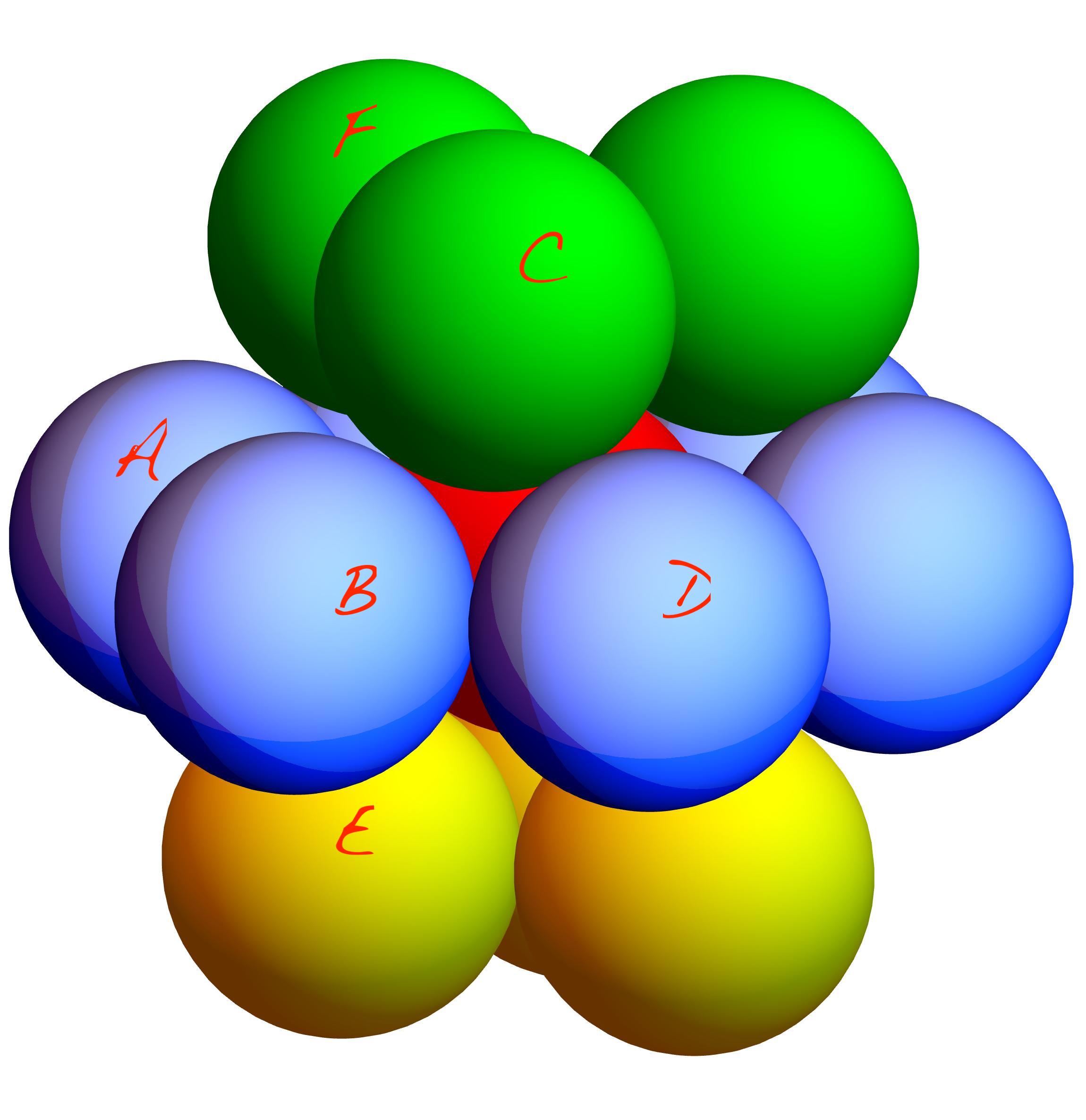}\hspace{2cm}
\includegraphics[scale=0.068]{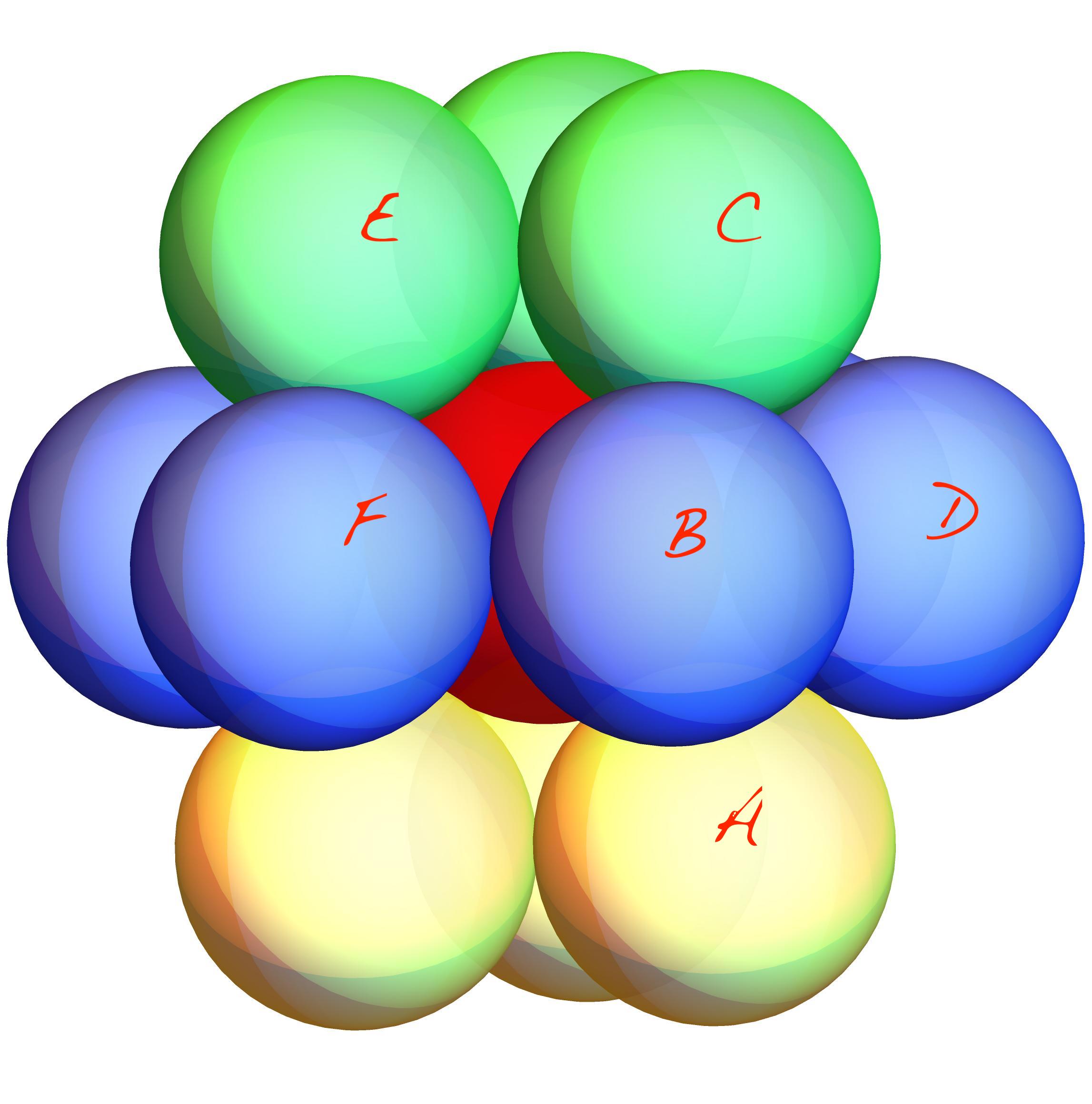}
\caption{FCC and HCP configurations}
\label{FCCandHCP}
\end{figure}

 The question is: can one roll the 12 balls over $\mathbb{B}$ in such a way that
they cease to touch (or \textit{kiss}) each other? (The
number of kissing points of the outer spheres is 24 for both configurations.)
When the answer is positive, we say that the configuration can be unlocked.

\vskip .2cm
The fact that the configuration FCC can be unlocked is well known. One of
unlocking moves is described in \cite{CS}. The situation with HCP cluster is
less known. In \cite{T} it is stated that the configuration HCP is rigid. But
in fact it is not the case, and the unlocking moves for both FCC and HCP are
presented in \cite{KKLS}.

\vskip .2cm
We present quite simple moves which unlock FCC and HCP.

\paragraph{FCC configuration.}The balls of the FCC configuration are centered at the vertices of the
cuboctahedron.

\begin{figure}[H]
\vspace{-.0cm}
\centering
\includegraphics[scale=0.46]{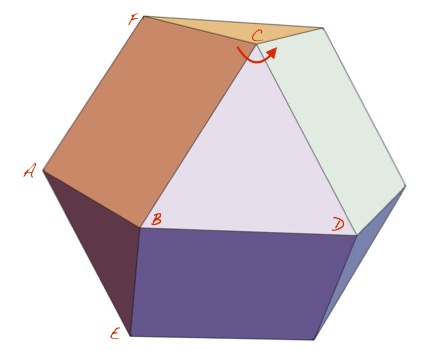}
\caption{FCC unlocking}
\label{FCCmove}
\end{figure}

\vskip .1cm
During the move the top three balls and the bottom three balls
remain fixed. Consider the three ``triangles of balls" sharing one ball with the top
triangle. (One of them -- the triangle $BCD$ -- is visible in Fig. \ref{FCCmove}.)
Now rotate each of these three
triangles with the same velocities around their top vertices, e.g. the triple
$C,B,D$ is rolled over the central ball, as a solid, with the ball $C$ fixed
(Fig. \ref{FCCmove}).

To see that the balls do not hit each other during the move, consider the ball $B$, which kisses four balls $A$, $C$, $D$ and $E$. The triple $BCD$ moves as a solid, so $B$ keeps kissing $C$ and $D$, and no conflict happens here. The ball $B$, being rotated around $C$, goes below the equator, while $A$, being rotated around $F$, goes above the equator, with the same equatorial projection of speeds, so the distance between them increases. The balls $C,B$ and $E$ lie on the same great circle, and both $C $ and $E$ do not move, so the distance $BE$ also increases.

\vskip .1cm
These considerations suffice due to the rotation symmetry (by $2\pi/3$
angle around the vertical axis) of the cuboctahedron.

\paragraph{HCP configuration.}The balls of the HCP configuration are centered at the vertices of the Johnson
solid called the triangular orthobicupola. They can be seen as three `rhombic'
configurations, one of them being $ABCD$ on the picture:

\begin{figure}[H]
\vspace{-.0cm}
\centering
\includegraphics[scale=0.46]{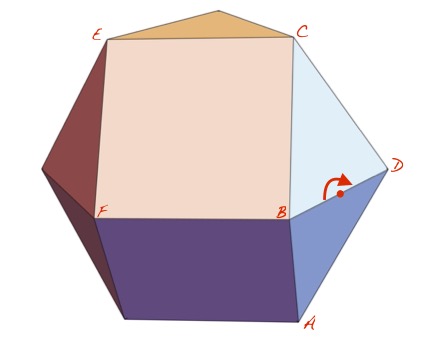}
\parbox{11.8cm}{
\caption{HCP unlocking} }
\label{HCPmove}
\end{figure}

\noindent Here $A,B,C,D$ are the centers of
the corresponding four balls. The move consists in rolling each of three
rhombi around their centers, with equal velocities. That is, one applies to
the four balls of each rhombus the rotation around an axis connecting the
center of the rhombus (the mid of the segment $BD$ for the rhombus
$ABCD$) and the origin.

\vskip .1cm
To see the absence of conflicts here, we observe that the balls $E$ and $C$ move downward, with the same horizontal projections of velocities, so their mutual distance increases (since their common distance to the north  pole increases).
The ball $B$ goes up, while the ball $F$ goes down, so they avoid each other.
The rest again follows from the same rotational $\mathbb{Z}_3$ symmetry of the orthobicupola.

\vskip .2cm
In both cases, FCC and HCP, after unlocking, the twelve balls can be moved apart and be positioned at the
vertices of the regular icosahedron. Then some free space appears between the unit
balls; their radii can be blown up to the value $r=\left(  \sqrt{\frac
{5+\sqrt{5}}{2}}-1\right)  ^{-1}\approx1.10851$, when they finally start to
kiss, making 30 kissing points. Yet this extra space is not large enough to
incorporate the 13-th unit ball, as was shown in \cite{SW}, settling the
famous discussion started by Newton and Gregory.

\section{Cylinders}

Our interest in the cylinder case started from the question of Kuperberg,
\cite{K}: one can easily put six unit cylinders around a central unit ball
$B:$

\begin{figure}[H]
\vspace{-.0cm}
\centering
\hspace{.8cm}\includegraphics[scale=0.4]{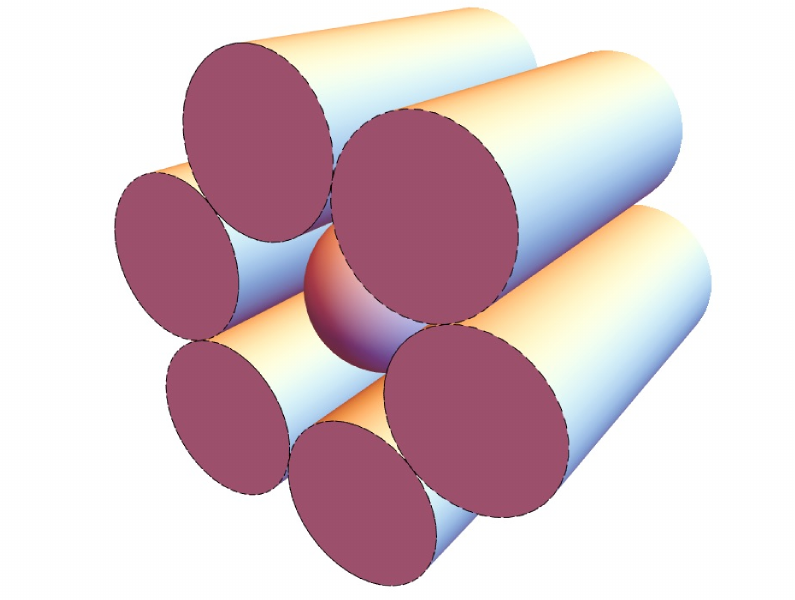}
\caption{Configuration $C_6$}
\label{SixParallelCylinders}
\end{figure}

The question is whether one can arrange seven such cylinders around $B$ in a
non-intersecting way. The question looks as an insult to intuition, but to
this day we do not have a rigorous answer to it. However one cannot
arrange eight unit cylinders \cite{BW}).

\vskip .1cm
The configuration shown looks quite solid; though it is not rigid --

\begin{figure}[H]
\vspace{-.2cm}
\centering
\hspace{.8cm}\includegraphics[scale=0.48]{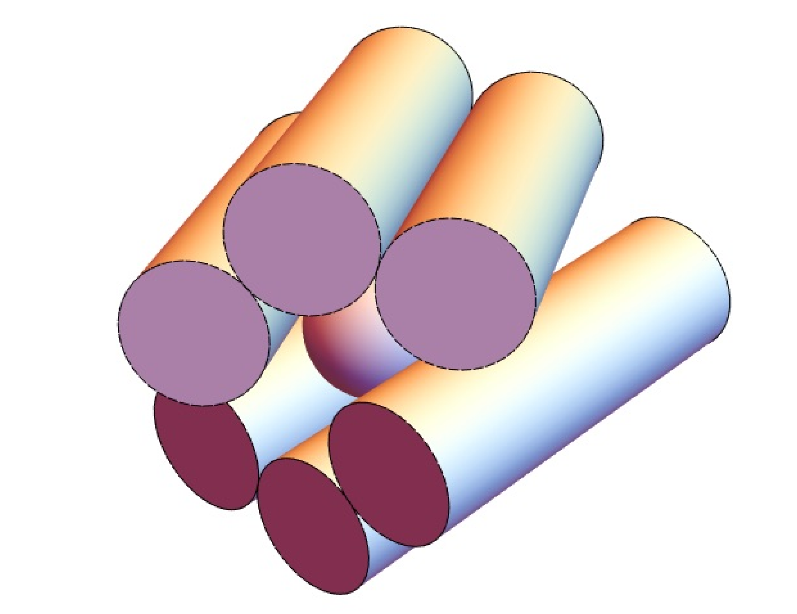}
\caption{Non-rigidity of $C_6$}
\label{NonRigidity}
\end{figure}

\vspace{-.2cm}
\noindent -- it looks pretty tight. So it came as a big surprise that there are
configurations of six unit cylinders around the unit ball, where the cylinders
do not touch each other. Such configurations were found by M. Firsching,
\cite{F}. In his example he was able to position six cylinders of the
radius $r=1.049659$ around the unit ball. This example was obtained by a
numerical exploration of the corresponding 18-dimensional configuration manifold.

\vskip .1cm
In our paper \cite{OS1} we found a way to unlock the configuration $C_6$ in a
symmetric (with respect to the dihedral group $\mathbb{D}_{3}$) manner, which enabled us to
improve the Firsching's $1.049659$ to the value
\begin{equation}
r_{\mathfrak{m}}=\frac{1}{8}\left(  3+\sqrt{33}\right)  \approx1.093070331...
\label{01}
\end{equation}
The corresponding $\mathbb{D}_{3}$-symmetric configuration $C_{\mathfrak{m}}$
is shown below:

\begin{figure}[H]
\vspace{-.4cm}
\centering
\hspace{.8cm}\includegraphics[scale=0.9]{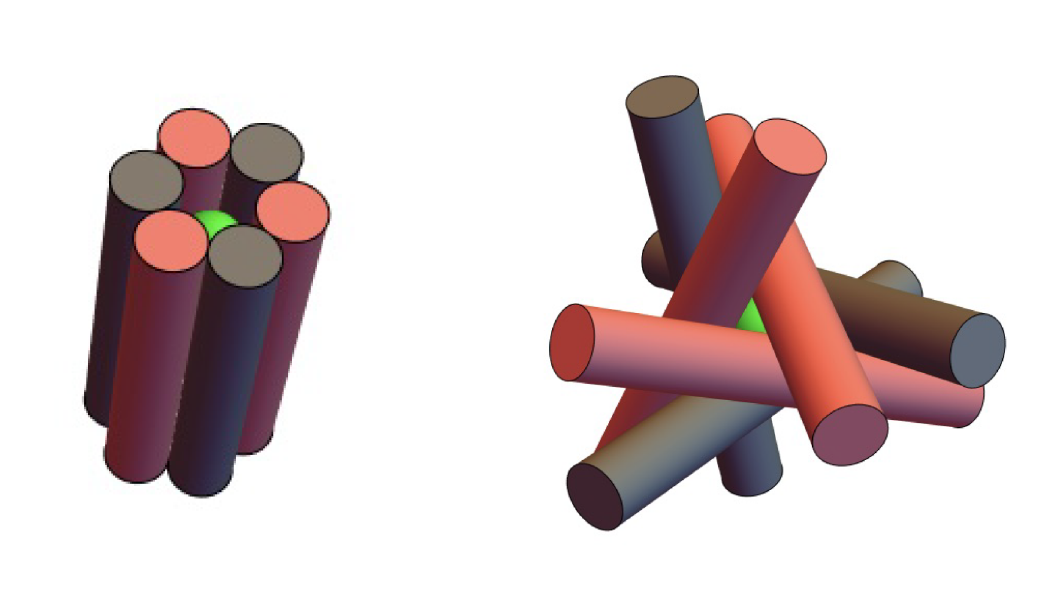}
\vspace{-.4cm}
\caption{The initial configuration $C_6$ and the record configuration
$C_{\mathfrak{m}}$}
\label{recordconf}
\end{figure}

The unlocking move looks as follows. Assume that the cylinders
on the left picture of Fig. \ref{recordconf} (the configuration $C_6$) point at the
North. We first describe a three-dimensional family of moves. The triple of pink cylinders go upward
by $\varphi,$ `horizontally' by $-\varkappa,$ and the three vectors
$\uparrow$ are rotated (around the axis joining the origin and the tangent point, counterclockwise if viewed from its tip) by $\delta,$ while the three remaining cylinders go
downward by $\varphi$, `horizontally' by $\varkappa,$ and the
three vectors $\uparrow$ are rotated, in the same way as the upper ones, by $\delta$. Now, optimising $\varkappa$ and $\delta$ for each $\varphi$,
we obtain a curve $\gamma$ in the moduli space of six cylinders. For each point on the trajectory $\gamma$ the configuration is $\mathbb{D}_{3}$-symmetric: its symmetries are 120$^{o}$ rotation around the axis North-South and 180$^{o}$ rotation around a perpendicular axis. The common radius of the cylinders grows as the configuration moves along $\gamma$ up to a certain point
which corresponds to the configuration $C_{\mathfrak{m}}$ and then decreases.
It turns out that the trajectory $\gamma$ is an {\bf algebraic}
curve (given by the relations (25)-(27) in \cite{OS1}), and this is the reason for the exact equality
$\left(  \ref{01}\right)  $. Y. Kallus made a movie
featuring our path $\gamma$, \cite{Ka}.

\vskip .1cm
It is interesting to
note that all the angles describing the configuration $C_{\mathfrak{m}}$ are
pure geodetic, in the sense of \cite{CRS}: an angle $\alpha$ is pure geodetic
if the square of its sine is rational.

\vskip .2cm
We believe that our configuration $C_{\mathfrak{m}}$ is the record
configuration, and in any configuration the radii of six equal
non-intersecting cylinders tangent to the unit ball $B$ are $\leq
r_{\mathfrak{m}}.$ But we have no proof of this maximality. We can prove,
however, a local version of this statement: if $C$ is a \textit{small proper}
perturbation of the configuration $C_{\mathfrak{m}},$ then the radii $r$ of
its cylinders satisfy: $r<r_{\mathfrak{m}}.$ Here \textit{proper }means that
$C$ is not a rotation of $C_{\mathfrak{m}},$ and \textit{small }means that the
six tangent points and six directions of the (equal) cylinders comprising $C$
are close to those of $C_{\mathfrak{m}}.$ The proof of this statement is a
content of our paper \cite{OS2}. The proof is quite
involved, since the function $r$ on our 18-dimensional configuration manifold
is not smooth. It turns out that the differentials of the distances between the tangent (to the unit central sphere) generatrices of the cylinders obey a convex linear dependence $\Lambda$ at the point $C_{\mathfrak{m}}$ on the curve
$\gamma$. Let $E$ be the linear
subspace of the tangent space on which all differentials vanish. Then we have proved that a sufficient condition for a local maximum is as follows: the {\bf same} linear combination $\Lambda$ of the {\it second}
differentials is negatively defined on $E$. We were able to check that this sufficient
condition holds at the point $C_{\mathfrak{m}}$ \cite{OS1}.

\vskip .2cm
{\bf Remark.} {\it Let $F$ denote the function $\max (F_1,\ldots,F_k)$ where
$F_1,\ldots,F_k$ are analytic functions on a vector space, and
$F_1(x_0)=\ldots =F_k(x_0)$ at some point $x_0$. The condition ``the restriction
of $F$ to any straight line passing through the point $x_0$ has a local maximum at
$x_0$'' is not sufficient to guarantee that the function $F$ has a local maximum at the point $x_0$, see example in \cite{OS2}, section 5. However we believe that if
the restriction of $F$ to any {analytic curve} passing through the point $x_0$ has a local maximum at
$x_0$ then the function $F$ does have a local maximum at the point $x_0$.

\vskip .1cm
Without the analyticity assumption on $F_{i}$-s the above conjecture does not
hold. We present an example of just one $C^{\infty}$ function $\Phi$ whose
restriction to any analytic path passing through $\left(  0,0\right)
\in\mathbb{R}^{2}$ has a local minimum at the origin, while the origin is not
a local minimum of $\Phi$. Consider the function
\[
\psi\left(  x\right)  =\left\{
\begin{tabular}
[c]{ll}
$\exp\left(  -\frac{1}{x}\right)  $ & if $x>0$\\
$0$ & if $x\leq0$
\end{tabular}
\ \ \right.  ,
\]
and let $\eta\left(  x\right)  $ be any $C^{\infty}$ function on
$\mathbb{R}^{1}$ with support on the segment $\left[  -\frac{1}{2},\frac{1}
{2}\right]  ,$ which is \noindent\textbf{i) }even, \textbf{ii) }decaying on
$[0,\infty),$ \textbf{iii) }$\eta\left(  0\right)  =1.$ We define the
$C^{\infty}$ function $\Phi$ on $\mathbb{R}^{2}$ by
\[
\Phi\left(  x,y\right)  =\left\{
\begin{tabular}
[c]{ll}
$\exp\left(  -\frac{1}{x^{2}}\right)  \eta\left(  \frac{y-\psi\left(
x\right)  }{\psi\left(  x\right)  }\right)  $ & for $x>0$\\
$0$ & for $x\leq0$
\end{tabular}
\ \ \right.  .
\]
The function $\Phi$ is indeed of $C^{\infty}$ class. The only problematic point is the origin $\left(  0,0\right)$. Any partial derivative (as $x\rightarrow 0^+,y\rightarrow 0$) of $\Phi$ belongs to the vector space (invariant with respect to the partial derivatives) of finite
linear combinations of functions $e^{-x^{-2}+a_1 x^{-1}}y^{a_2}x^{-a_3}
\eta^{(a_4)}\left(  ye^{x^{-1}}-1\right)$ for some $a_1,a_2,a_3,a_4\in
\mathbb{Z}_{\geq 0}$ and hence vanishes.

Let us show that for any analytic path $\gamma:\left[  0,1\right]
\rightarrow\mathbb{R}^{2}$, $\gamma\left(  0\right)  =\left(  0,0\right)  ,$
the function $\Phi\left(  \gamma\left(  t\right)  \right)  $ vanishes on some
segment $0\leq t<u_{\gamma},$ so $0$ is a local minimum of the function
$\Phi\left(  \gamma\left(  t\right)  \right)  $ on $\left[  0,1\right]  ,$
while the point $\left(  0,0\right)  \in\mathbb{R}^{2}$ is evidently not a
local minimum of the function $\Phi.$ To see this, note that the support of
the function $\Phi$ lies inside the `beak' $\beta=\left\{  \left(  x,y\right)
:x\geq0,\frac{1}{2}\psi\left(  x\right)  \leq y\leq2\psi\left(  x\right)
\right\}  \subset\mathbb{R}^{2}.$ Let $\gamma\left(  t\right)  =\left(
x\left(  t\right)  ,y\left(  t\right)  \right)  $ be an analytic path, defined
by two analytic functions $x\left(  t\right)  ,y\left(  t\right)  .$ Then
$x\left(  t\right)  =a^{\prime}t^{k^{\prime}}+O\left(  t^{k^{\prime}
+1}\right)  ,$ $y\left(  t\right)  =a^{\prime\prime}t^{k^{^{\prime\prime}}
}+O\left(  t^{k^{^{\prime\prime}}+1}\right)  $ for some real $a^{\prime
},a^{\prime\prime}$ and integer $k^{\prime},k^{\prime\prime},$ and therefore
the path $\gamma$ in the vicinity of the origin is a graph of a function
$y_{\gamma}\left(  x\right)  =bx^{k^{^{\prime\prime}}/k^{\prime}}+o\left(
x^{k^{^{\prime\prime}}/k^{\prime}}\right)  $ for some real $b$. (Without loss
of generality we can assume that $k^{\prime}>0.)$ Therefore there exists a
value $\varepsilon_{\gamma}>0$ such that the graph of the function $y_{\gamma
}$ does not intersect the beak $\beta$ for $0<x<\varepsilon_{\gamma},$ because
the function $\exp\left(  -\frac{1}{x}\right)  $ is smaller than any power of
$x$ in the appropriate vicinity of $0.$ Therefore the function $\Phi$ vanishes
on that piece of $\gamma.$

Apparently, to probe some $C^{\infty}$ function, one needs all $C^{\infty}$
paths, and not just analytic paths.

}

\vskip .2cm
The local maximality of the configuration $C_{\mathfrak{m}}$ implies that it
is rigid, i.e. cannot be unlocked. It seems to play the role of the
icosahedral configuration of 12 kissing balls above, while the six unit
cylinders can be rolled away from each other and create some free space
between them on the sphere. Yet whether this space is sufficient for the
seventh unit cylinder to be squeezed in ( $\equiv$ Kuperberg question) is
unknown, in contrast to the 13 unit balls problem. The thesis of O. Yardimci,
\cite{Y}, contains the theorem, proven together with A. Bezdek, that one
cannot put 7 cylinders of the Firsching radius $1.049659$ in contact with the
unit central ball, in a non-intersecting way.

\vskip .2cm
In addition to configurations on Figs. \ref{SixParallelCylinders} and
\ref{NonRigidity}, W. Kuperberg pointed out yet
another configuration of six unit cylinders around the unit ball:

\begin{figure}[H]
\vspace{-.4cm}
\centering
\hspace{.8cm}\includegraphics[scale=0.7]{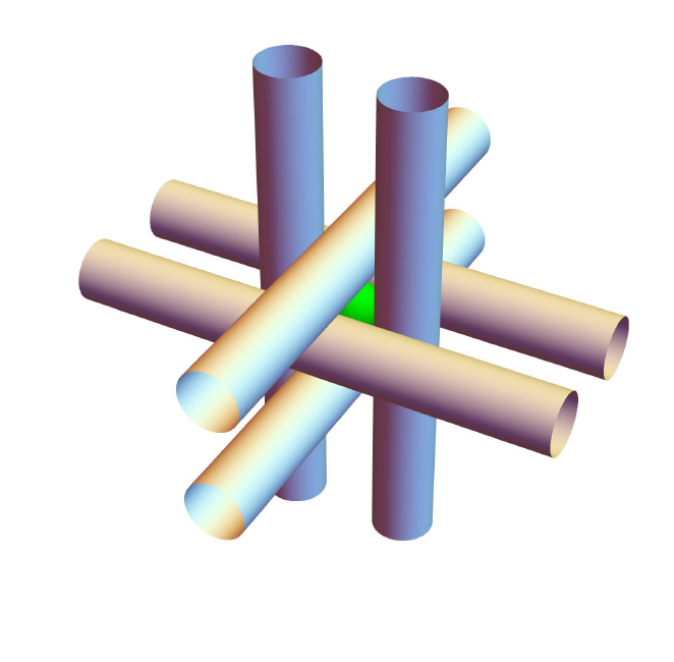}
\vspace{-.8cm}
\caption{Configuration $O_6$}
\end{figure}

\vspace{-.2cm}
It looks rigid, and in our paper
\cite{OS3} we have proven this rigidity: for every \textit{small proper}
perturbation of the configuration $O_{6}$, the radii $r$ of its cylinders
satisfy: $r<1.$ In particular, the configuration $O_{6}$ also cannot be
unlocked. The proof turned out to be even more involved.
The differentials of the distances between the tangent (to the unit central sphere) generatrices of the cylinders obey three convex linear dependencies
$\Lambda_1$, $\Lambda_2$ and $\Lambda_3$. Let again $E$ be the linear
subspace of the tangent space on which all differentials vanish and denote
by $q_1,q_2,q_3$ the same linear combinations $\Lambda_1$, $\Lambda_2$ and $\Lambda_3$ of the second derivatives. We have checked that the system $q_1(x)>0,q_2(x)>0,q_3(x)>0$ of inequalities has no solution on $E$ and
proved that this is a sufficient condition for a local maximum \cite{OS3,OS5}.
Interestingly, no convex linear combination of the forms $q_1,q_2,q_3$ is negatively defined on $E$ so we could not just apply the technics like Sylvester criterion. We think that it is a beautiful program -- to develop the theory (examples, criteria, classification) of tuples
$q_1,\ldots ,q_j$ of quadratic forms on a vector space for which the system $q_1(x)>0,\ldots ,q_j(x)>0$ of inequalities has no solution.

\vskip .1cm
Sometimes this configuration $O_{6}$ is called `octahedral' because, probably,
the tangency points lie at the vertices of the regular octahedron. We present
an interpretation of the configuration $O_{6}$ which relates it to the
configuration of rotated edges of the regular tetrahedron. This
interpretation, which shows that it rather deserves to be called the
tetrahedral configuration is as follows. Consider the configuration of the
tangent to the unit sphere lines which are continuations of the edges of the
regular tetrahedron. The points of the sphere at which tangent lines pass are
the edge middles of the tetrahedron:

\begin{figure}[H]
\vspace{-.2cm}
\centering
\hspace{.8cm}\includegraphics[scale=0.7]{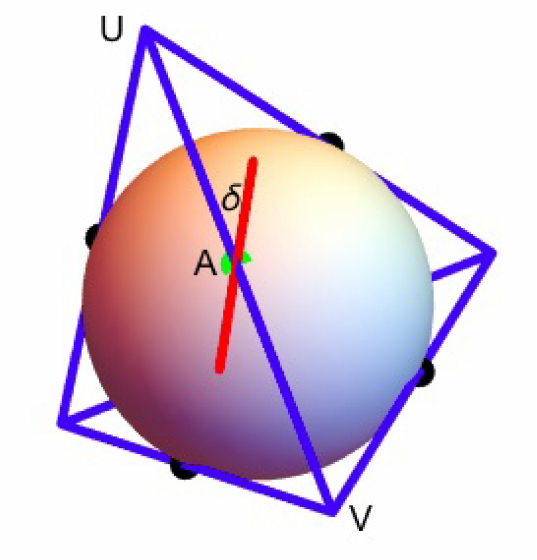}
\caption{$\delta$-process}
\end{figure}

\vspace{-.2cm}
Then each edge is rotated around the diameter of the unit sphere, passing
through the middle of the edge, by an angle $\delta$. On our picture the point
$A$ (in green) is the middle of the edge $UV$. The line, passing through the
point A and rotated by the angle $\delta$, is shown in red. The lines passing
through other five middles of edges are rotated by the same angle $\delta,$ in
accordance with the group $\mathbb{A}_{4}$ of the proper symmetries of the tetrahedron. This motion we call {\it $\delta$-process}.

\vskip .1cm
Let is replace each rotated line by a cylinder of radius $r,$ tangent to the
sphere, in such a way that the line is its generatrix, and (some of) the
cylinders are kissing, which uniquely defines $r$ as a function $r\left(
\delta\right)  $ of the angle $\delta$. Then for $\delta=0,$ $\pi/2$ the
radius $r=0,$ while it is maximal at $\delta=\pi/4,$ with $r\left(
\pi/4\right)  =1.$ This is precisely the configuration $O_{6}.$

\vskip .1cm
Following this interpretation we have introduced in \cite{OS4} configurations
of tangent cylinders for two remaining pairs of dual Platonic bodies, that is,
for the octahedron/cube $\left(  \mathcal{O}/\mathcal{C}\right)  $ and
icosahedron/dodecahedron $\left(  \mathcal{I}/\mathcal{D}\right)  $; for the
pair tetrahedron/tetrahedron this is precisely the configuration $O_{6}$. The
number of cylinders in our configurations is equal to the number of edges of
either of Platonic bodies in the pair, that is, twelve for the pair
octahedron/cube and thirty for the pair icosahedron/dodecahedron. The
corresponding radii $r_{\mathcal{O}/\mathcal{C}},\ r_{\mathcal{I}/\mathcal{D}
}$ of the cylinders are obtained as a result of the maximisation of the functions
$r_{\mathcal{O}/\mathcal{C}}\left(  \delta\right)  ,\ r_{\mathcal{I}
/\mathcal{D}}\left(  \delta\right)  $ over the rotation angle $\delta.$

\vskip .2cm
\noindent{\bf Pair $\mathcal{O}/\mathcal{C}$.}
For the pair $\mathcal{O}/\mathcal{C}$ the optimal configuration is the following:

\begin{figure}[H]
\vspace{-.4cm}
\centering
\includegraphics[scale=0.24]{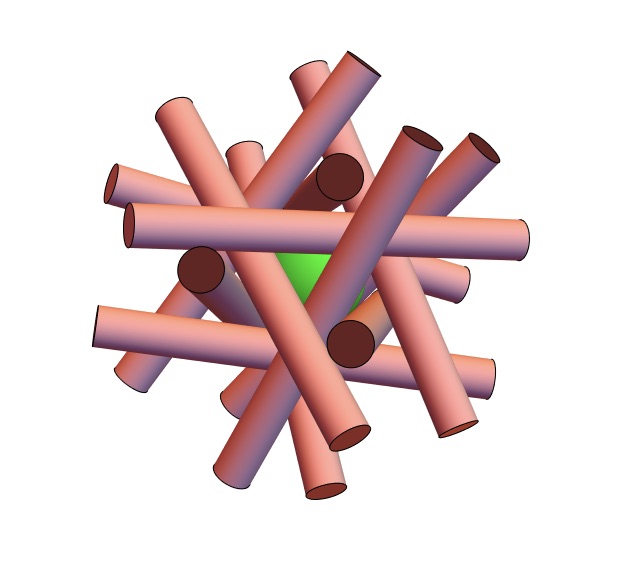}
\vspace{-.4cm}
\caption{$\mathcal{O}/\mathcal{C}$ maximal configuration of cylinders, view from a vertex of the cube}
\label{OC-1}
\end{figure}

\vspace{-.3cm}
\noindent Each cylinder touches four other cylinders. The corresponding value
$\delta_{\mathcal{O}/\mathcal{C}}$ satisfies
$
\tan(\delta_{\mathcal{O}/\mathcal{C}})=\frac{3^{1/4}}{\sqrt{2}},
$
(approximately, $\delta_{\mathcal{O}/\mathcal{C}}\simeq0.23856\pi
\simeq0.74946$). The corresponding radius of touching cylinders is
\[
r_{\mathcal{O}/\mathcal{C}}=\frac{\sqrt{3}-1}{1+2\sqrt{2}-\sqrt{3}}
\approx0.3492\ .
\]

\vskip .2cm
\noindent{\bf Pair  $\mathcal{I}/\mathcal{D}$.} For the pair $\mathcal{I}/\mathcal{D}$ the optimal configuration is the following:

\begin{figure}[H]
\centering
\includegraphics[scale=0.3]{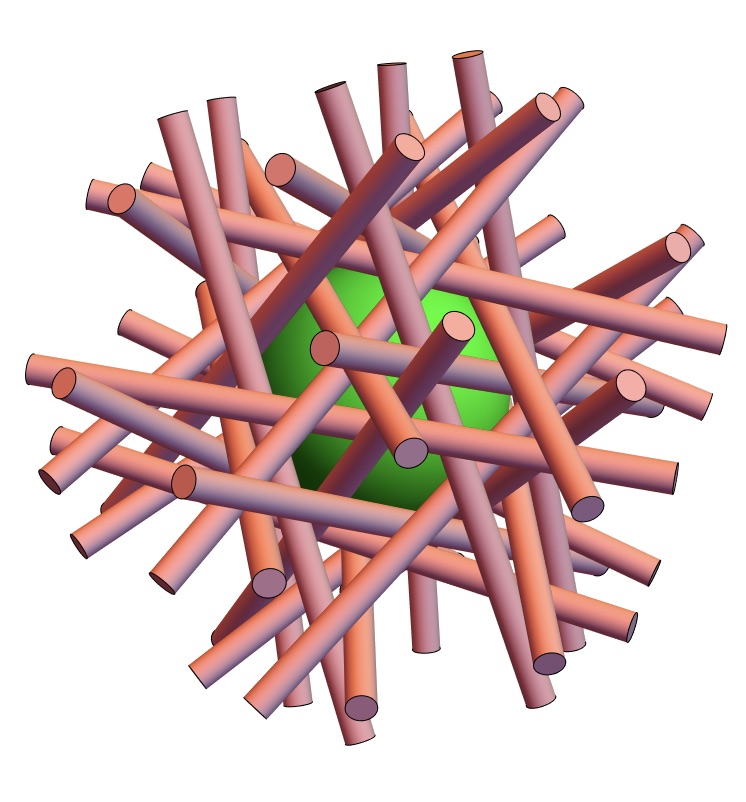}
\label{MaxID-ConfigViewFrom3-axis}
\caption{Maximal $\mathcal{I}/\mathcal{D}$ configuration, view from the tip of a 3-fold axis}
\end{figure}

\vspace{-.4cm}
\noindent
Here each cylinder touches eight other cylinders. The optimal angle
$\delta_{\mathcal{I}/\mathcal{D}}$ is given by
\[
\delta_{\mathcal{I}/\mathcal{D}}=\arctan\left(  \sqrt{t_{0}}\right)  \ ,
\]
where $t_{0}\simeq0.694356$ and is a root of the polynomial
\[
5t^{6}-80t^{5}+190t^{3}-4t^{2}-84t+9\ .
\]
Approximately,
\[\delta_{\mathcal{I}/\mathcal{D}}\simeq0.694707\ .\]
The
corresponding radius of cylinders is approximately
\[ r_{\mathcal{I}
/\mathcal{D}}\simeq0.115558\ .\]
Another image of this configuration:

\begin{figure}[H]
\centering
\includegraphics[scale=0.32]{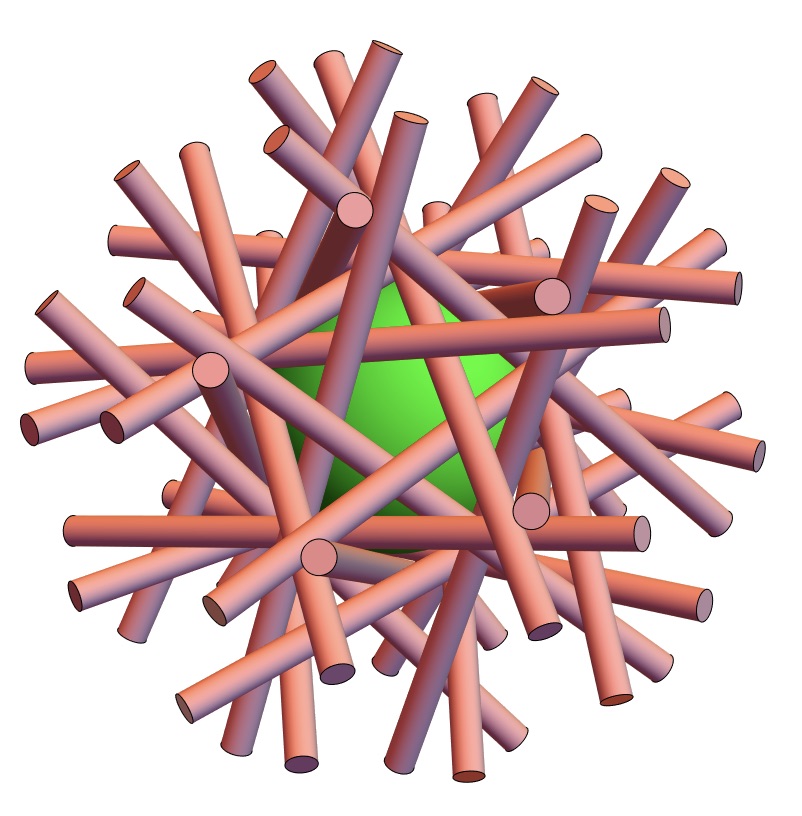}
\label{MaxID-ConfigViewFrom5-axis}
\caption{Maximal $\mathcal{I}/\mathcal{D}$ configuration, view from the tip of a 5-fold axis}
\label{zzve}
\end{figure}

\noindent is shown on the cover.

\vskip .2cm
We conjecture that both record configurations -- $\mathcal{O}/\mathcal{C}$ and
$\mathcal{I}/\mathcal{D}$ -- are rigid. This probably can be checked by using
our machinery developed in \cite{OS2} and \cite{OS3}, but the computations are
too lengthy, and we have not attempted to perform them.

\vskip .2cm
Another interesting phenomena happens during the $\delta$-process at the angle
values when the function $r_{\mathcal{I}/\mathcal{D}}\left(  \delta\right)  $
vanishes. That means that some of the 30 lines are intersecting. There are
three values of $\delta$ at which this happens. The patterns emerging are: ten
linked triangles, five linked tetrahedra and six linked pentagonal stars.

\begin{figure}[H]
\vspace{-.4cm} \centering
\hspace{-.6cm}
\includegraphics[scale=0.28]{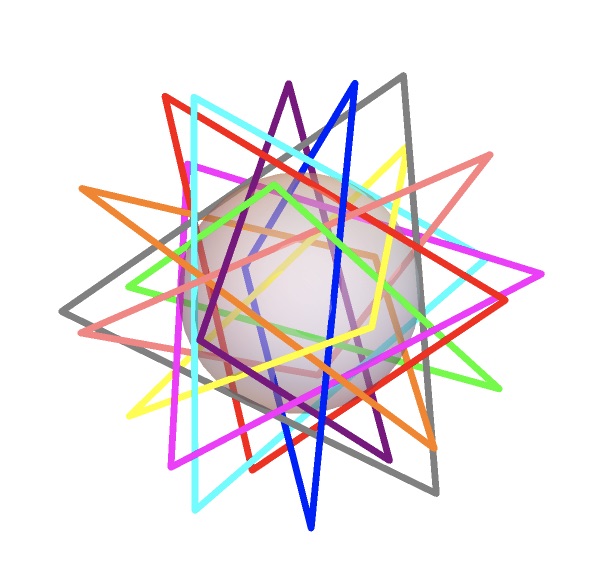}\hspace{-.7cm}
\includegraphics[scale=0.36]{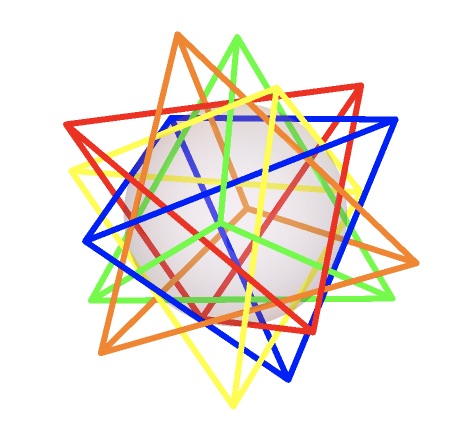}\hspace{-.7cm}
\includegraphics[scale=0.3]{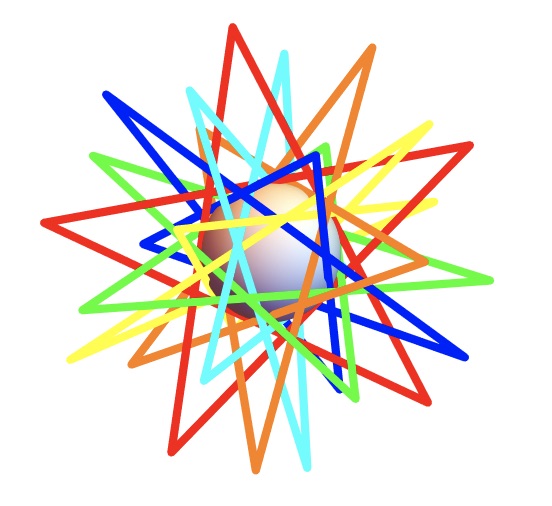}
\caption{$\mathcal{I}/\mathcal{D}$ minima}
\end{figure}

\end{document}